\documentclass[12pt]{amsart}
\usepackage{amsmath, amsthm, amscd, amssymb, amsfonts, amsxtra, amssymb, latexsym}
\usepackage{enumerate}
\usepackage{verbatim}
\usepackage{array}
\usepackage[utf8]{inputenc}
\usepackage{epsfig}
\usepackage{graphicx,tikz}
\usepackage{shadow}
\usepackage{rotating}
\usepackage{wrapfig,sidecap}
\usepackage{multicol}
\usepackage{pst-node}
\usepackage{pstricks}

\usepackage{hyperref}
\hypersetup{colorlinks = true,	allcolors  = blue}

\def\fq{\mathbb{F}_{q}}
\def\fqs{\mathbb{F}_{q^2}}

\def\cC{\mathcal{C}}

\def\cP{\mathbb{P}}

\def\cL{\mathcal{L}}

\def\N{\mathbb{N}}

\newcommand{\con}{\operatorname{Con}}

\newcommand{\rr}{\mathcal{L}}
\newcommand{\sop}{\operatorname{Supp}}

\newcommand{\ff}{\operatorname{\mathbb{F}}}

\def\b{\beta}

\def\ff{\mathbb{F}}

\DeclareMathOperator\supp{supp}
\DeclareMathOperator\diff{Diff}

\usepackage{color}

\newcommand{\sk}{\smallskip}

\newtheorem{theorem}{Theorem}[section]

\newtheorem{corollary}[theorem]{Corollary}
\newtheorem{proposition}[theorem]{Proposition}

\theoremstyle{definition}
\newtheorem{definition}[theorem]{Definition}

\newtheorem{remark}[theorem]{Remark}

\begin{document} \sloppy
	\numberwithin{equation}{section}
	\title[Lifting canonical divisors in towers]{Lifting canonical divisors of the form $2G-D$ in towers of function fields} 
	\author[M.\@ Chara, R.\@ Podestá, L.\@ Quoos, R.\@ Toledano]{María Chara, Ricardo Podestá, Luciane Quoos, Ricardo Toledano}
	\dedicatory{\today}
	\address{María Chara, Researcher of {\sc Pedeciba} at Facultad de Ciencias, Universidad de la República, Iguá 4225, (11400) Montevideo, Uruguay. 
	\newline {\it E-mail: mchara@cmat.edu.uy}}
	
	\address{Ricardo Podest\'a, FaMAF -- CIEM (CONICET), Universidad Nacional de C\'ordoba, \newline
	Av.\@ Medina Allende 2144, Ciudad Universitaria, (5000) C\'ordoba, Argentina. 
	\newline {\it E-mail: ricardo.podesta@unc.edu.ar}}
	
	\address{Luciane Quoos, Universidade Federal do Rio de Janeiro, Centro de Tecnologia, Cidade Universitária, Av.\@ Athos da Silveira Ramos 149, Ilha do Fundão, CEP~21.941-909, Brazil. \newline {\it E-mail: luciane@im.ufrj.br}}
	
	\address{Ricardo Toledano, Facultad de Ingeniería Química, Universidad Nacional del Litoral, Santiago del Estero 2829, (3000) Santa Fe, Argentina.\newline {\it E-mail: ridatole@gmail.com}}
	\dedicatory{\today}

\begin{abstract}
	We present a method for constructing, given a canonical divisor in a function field of the form $2G-D$, another canonical divisor of the form $2G'-D'$ in an extension of the function field that satisfies certain ramification conditions. This method applies in particular to certain recursive towers that yield sequences of canonical divisors of the form $2G_i-D_i$ for each index  $i\geq 0$. Using these divisors, we construct sequences of iso-dual AG-codes and establish the existence of asymptotically good sequences among them. Beyond its simplicity, a key advantage of our approach is its flexibility, which allows it to be applied to a wide range of recursive towers of function fields, thereby providing a more explicit way to construct asymptotically good sequences of iso-dual AG-codes.
\end{abstract}

\maketitle

\section{Introduction} \label{sec:intro}
Let $\fq$ be a finite field with $q$ elements. A \textit{linear code} $\cC$ over $\fq$ is an $\fq$-linear subspace of $\fq^n$ for $ n\geq 1$. Associated to a linear code $\cC \subseteq \fq^n$ there are three important parameters: its \textit{length} $n$, its \textit{dimension} $k$ as a vector space over $\fq$ and its \textit{minimum} (Hamming) \textit{distance} $d$. 
To emphasize these parameters, it is usual to speak of a linear code $\cC$ as an $[n, k, d]$-code over $\fq$.
The importance of these parameters is that they tell us how good a linear code is for detection and correction of errors. 
In the 1980's, using tools coming from algebraic geometry, Goppa constructed linear codes using algebraic curves defined over a finite field, see \cite{G1981} and \cite{G1982}. The codes constructed by Goppa are usually called algebraic geometry (AG) codes and have played an important role in Coding theory, mainly because Goppa's ideas were used by Tsfasman, Vlăduţ and Zink  in \cite{TVZ1982} to show that there are sequences of linear codes whose parameters are asymptotically better than the Gilbert–Varshamov bound. This certainly was a groundbreaking result at that time. Another important result is that every linear code can be realized as an algebraic geometry code \cite{PSV1991}.

The dual code $\cC^\perp$ of a linear code $\cC$ is the orthogonal complement of $\cC$ in $\fq^n$ with respect to the standard (Euclidean) inner product of $\fq^n$. A code is said to be self-dual if $\cC = \cC^\perp$.
Self-dual AG-codes began to be systematically investigated in \cite{Stich1988}, with generalizations and new constructions proposed in, for example, \cite{DM03} and \cite{S2021}. To ensure that an AG code $C_\rr(D,G)$ be self-dual  strong conditions must be imposed as it is shown in Corollary 2.2.11 of \cite{Stichbook09}: a Weil differential $\eta$ has to be found such that $(\eta)=2G-D$, so that $2G-D$ has to be a canonical divisor, and  the local component $\eta_P(1)=1$ for each $P\in \supp D$. It is worth to recall that these are not necessary conditions for self-duality of AG-codes because there are examples of self-dual AG-codes $C_\rr(D,G)$ such that $2G-D$ is not a canonical divisor (see Example 4.1.2 of \cite{Stich1988} and also \cite{CPQT2024corr}). All of this shows, in particular, that canonical divisors of the form $2G-D$ are important objects of study.

A slightly more general and less restrictive situation  is achieved with the so called iso-dual codes. More precisely (see Definitions \ref{defequiv} and \ref{defisodual} in the  next section) a linear code is called iso-dual if it is equivalent to its dual. Iso-dual codes can be naturally defined not only over finite fields but also over finite rings. For example, they were used in \cite{BGAH00} to study the construction of iso-dual lattices, an interesting class of lattices introduced by Conway and Sloane in \cite{CS94}.

The main goal of this article is to provide a general and relatively simple method to construct sequences of canonical divisors of the form $2G-D$ in each step of a tower of function fields assuming some suitable conditions on the ramification behavior. As a direct consequence of our method we are able to obtain asymptotically good sequences of iso-dual AG-codes using certain recursive towers of function fields.  This is an important difference with the case of self-dual AG codes: the existence of good sequences of self-dual AG-codes have been proved in \cite{Stich06} and \cite{BS2019} but by using  the Galois closure of recursive towers of function fields which, as far as we know, are not recursively defined. Another advantage of our method is that it is not only simpler than the ones used in \cite{Stich06} and \cite{BS2019} but it is more flexible, in the sense that it can be used in a greater variety of towers of function fields including  those ones recursively defined which are the most important ones in the applications. In this direction we use our method to show in Proposition~\ref{prop: TVZ} the existence of a sequence of iso-dual AG-codes over non-prime fields of quadratic cardinality that attains the Tsfasman–Vlăduţ–Zink lower bound in a relatively simple way and using an optimal recursive tower defined in \cite{GS96}. We also use our method to show, in Proposition~\ref{prop: goodisodualF8}, the existence of good sequences of iso-dual AG-codes over finite fields of cubic cardinality and even characteristic by using the recursive tower defined in \cite{BGS05}. In this case the lower bound we obtain for the relative minimum distance of our good sequence of iso-dual AG codes is the same lower bound for the relative minimum distance of the good sequences of self-dual AG-codes given in \cite{BS2019}.

The organization of this paper is as follows. In Section \ref{sec:preliminaries}, we begin by giving the basic definitions and properties we will use from the theory of function fields as well as the definition of Weil differentials and canonical divisors. We also define the concept of lifting a divisor in an extension of function fields. We then recall the construction of AG-codes using function fields over a finite field and provide the definition of iso-dual codes. In Section \ref{sec:lifting} we recall what a tower of function field is and we recall in Theorem \ref{teolevantado} the main result of our previous work \cite{CPQT2024} (see also \cite{CPQT2024corr} for the correct version of the main result of \cite{CPQT2024}) giving the conditions under which it is possible to lift a canonical divisor to a finite extension of  function fields. We give conditions in Theorem \ref{liftingcandivtower} to  lift  a canonical divisor of the form $2G-D$ from the bottom field of a tower of function fields to any step of the tower. In Section \ref{sec: goodsequences} we recall what a good sequence of linear codes is  and we apply our method of lifting canonical divisors of the form $2G-D$ in a tower to prove the existence of good sequences of iso-dual codes. In Theorem \ref{isontowers}, we provide the conditions needed to ensure the existence of good sequences of iso-dual AG-codes using our previous results of Section \ref{sec:lifting} on the lifting of canonical divisors along a tower and we give in Corollary \ref{isoratexamples} our method to obtain good sequences of iso-dual AG-codes starting from a rational function field, which is the most useful in the applications to Coding theory.   We then use this method given in Corollary \ref{isoratexamples} in  well-known recursive towers of function fields: for example in Proposition \ref{prop: TVZ} the existence of a good sequence of iso-dual AG-codes over a finite field of quadratic cardinality that attains the well-known TVZ bound is proved by using Corollary \ref{isoratexamples} in an optimal recursive tower studied in \cite{GS96}. We  keep using our Corollary \ref{isoratexamples} in another recursive tower of function fields studied in \cite{BGS05} to show, in Proposition \ref{prop: goodisodualF8}, the existence of good sequences of iso-dual AG-codes over finite fields of cubic cardinality and even characteristic. We conclude this work with Section \ref{succliftings}, where we show that if a canonical divisor of the form $2G-D$, defined at the bottom field $F_0$ of a tower $\mathcal{T} = \{F_i\}_{i=0}^\infty$, is lifted to a canonical divisor of the $i$-th step $F_i$ of $\mathcal{T}$ directly, then this canonical divisor coincides with the one obtained by  lifting $2D-G$ step by step in the tower until reaching the $i$-th step $F_i$ of $\mathcal{T}$.

\section{Preliminaries} \label{sec:preliminaries}
Here we recall some basic definitions and facts of function fields, divisors, AG-codes and iso-dual AG-codes following the presentation given in \cite{Stichbook09}.

\subsection*{Algebraic function fields} 
Let $\fq$ be a finite field with $q$ elements and let $F$ be an algebraic function field over $\fq$. All function field extensions $F'/F$ considered in this work are always finite, separable and geometric. The last condition means that $\ff_q$ is the full constant field of both $F'$ and $F$.  

We denote by $\mathbb{P}_F$ the set of places in $F$ and by $\nu_{P}$ the discrete valuation of $F$ associated to the place $P\in \cP_F$. Associated to each valuation $\nu_P$ we have the discrete valuation ring $\mathcal{O}_P=\{x\in F:\nu_P(x)\geq 0\}$. 

By $\mathrm{Div}(F)$ we denote the free abelian group generated by the places in $F$. The elements of $\mathrm{Div}(F)$ are called {\it divisors} of $F$ and are denoted by $D=\sum_{P \in \mathbb{P}_F} a_P P$, where $a_P=0$ for almost all $P$. The support $\sop(D)$ of a divisor $D$ is the set of places such that $a_P\neq 0$. 

Any discrete valuation $\nu_P$ can be extended to the group $\mathrm{Div}(F)$ by defining $\nu_P(D)=a_P$ if $D=\sum_{P \in \mathbb{P}_F} a_P P$.

Two divisors $D_1$ and $D_2$ of $F$ are called {\it disjoint divisors} if $\sop(D_1)\cap \sop(D_2)= \varnothing$. For a function $z \in F$ we let 
	$$(z)^F, \quad (z)^F_\infty, \quad \text{and} \quad (z)^F_0$$ 
stand for the principal, the pole and the zero divisors of the function $z$ in $F$, respectively. If needed we will write $(z)^F$ to make clear that we are considering the principal divisor of $z$ in $F$. Two divisors $A,B \in \mathrm{Div}(F)$ are equivalent, denoted $A \sim B$, if they differ by a principal divisor, i.e.\@ $B=A+(z)^F$ for some $z \in F$.

Let us consider now a function field extension $F'/F$. A place $Q$ of $F'$ is {\it above} or {\it lying over} a place $P$ of $F$ if $Q\cap F=P$. In this case we also say that $P$ is {\it under} or {\it lying under} $Q$ and we write $Q|P$ when $Q\cap F=P$. A set of places $\{Q_1,\ldots,Q_r\}$ of $F'$ is called the {\it overset} of a set of places $\{P_1,\ldots,P_n\}$ of $F$ if 
\[\{Q_1,\ldots,Q_r\}\cap F=\{P_1,\ldots,P_n\}.\]
Let $D\in \mathrm{Div}(F)$ and $D'\in \mathrm{Div}(F')$. We say that $D'$ is the {\it lifting} of $D$ to $F'$ if $\sop(D')$ is the overset of $\sop(D)$.

Let $Q$ be a place of $F'$ lying above a place $P$ of $F$. The symbols $e(Q|P)$, $f(Q|P)$ and $d(Q|P)$ denote the \textit{ramification index}, the \textit{inertia degree} and the \textit{different exponent} of $Q$ over $P$, respectively.

Let $P$ be a place of $F$. It is said that $P$ \textit{splits completely} in $F'$ if there are exactly $[F':F]$ places of $F'$ lying over $P$. This is equivalent to have that $e(Q|P)=f(Q|P)=1$ for every place $Q$ of $F'$ lying over $P$. The opposite situation occurs when there is only one place $Q$ of $F'$ lying over $P$. In this case, if $e(Q|P)=[F':F]$ (hence $f(Q|P)=1$) it is said that $P$ is \textit{totally ramified} in $F'$ and if $f(Q|P)=[F':F]$ (hence $e(Q|P)=1$) it is said that $P$ is \textit{inert} in $F'$.

The \textit{conorm} of a place $P$ of $F$ is the divisor of $F'$ defined as
	$$ \con_{F' /F} (P) = \sum_{Q \in F', Q|P} e(Q|P) \, Q.$$
For a divisor $D\in \mathrm{Div}(F)$ the conorm map defines a group homomorphism from $\mathrm{Div}(F)$ to $\mathrm{Div}(F')$ 
when extended by linearity; 
that is, 
\begin{equation} \label{eq:conorm}
	\con_{F' /F} \Big( \sum_i n_i P_i \Big) = \sum_i n_i \, \con_{F' /F} (P_i).
\end{equation}

The \textit{different divisor} of the extension $F'/F$ is the divisor of $F'$ defined as 
\begin{equation} \label{eq:diff}
	\diff(F'/F) =\sum_{P \in \mathbb{P}_F} \sum_{Q|P} d(Q|P) \, Q,
\end{equation}
where $d(Q|P)$ denotes the different exponent of $Q$ over $P$.
The different divisor is well defined because $d(Q|P)$ is a non-negative integer which is zero  for almost all places.

An important formula relates the genera $g$ and $g'$ of $F$ and $F'$ respectively. This formula is usually called 
the Riemann-Hurwitz genus formula and it asserts that 
\begin{equation} \label{eq: genus fla}
	2g'-2 = [F':F](2g-2) + \deg(\diff(F'/F)).
\end{equation}
 
We recall now the definition of Weil divisors and canonical divisors of $F$. An {\it adele} of $F$ is a map
 \[\alpha:\mathbb{P}_F\rightarrow F,\]
such that $\alpha(P)\in \mathcal{O}_P$ for almost all $P\in \mathbb{P}_F$. We can identify each element $x\in F$ with an adele of $F$ by defining $\alpha(P)=x$ for all $P\in \mathbb{P}_F$. The collection of all adeles of $F$ is denoted as $\mathcal{A}_F$ and it is  a vector space over $\fq$ and also a ring. A discrete valuation $\nu_P$ of $F$ can be extended to $\mathcal{A}_F$ by defining $\nu_P(\alpha)=\nu_P\bigl(\alpha(P)\bigr)$.

A {\it Weil differential} of $F$ is an $\fq$-linear map 
\[\eta: \mathcal{A}_F\rightarrow \fq,\] 
such that $\eta(\alpha)=0$ for all $\alpha\in \mathcal{A}_F(D)+F$ where $\mathcal{A}_F(D)$ is the collection of all adeles $\alpha\in \mathcal{A}_F$ such that $\nu_P(\alpha)\geq -\nu_P(D)$ for all $P\in\mathbb{P}_F$. The collection of all Weil differentials is denoted as $\Omega_F$. For each $\eta\in \Omega_F$ there is the unique divisor of $F$, denoted as $(\eta)$, such that $\eta$ vanishes on $\mathcal{A}_F((\eta))+F$ and if $\eta$ vanishes on $\mathcal{A}_F(D)+F$ for some $D\in \mathrm{Div}(F)$ then $D\leq (\eta)$. 

A divisor $D$ of $F$ is called {\it canonical} if $D=(\eta)$ for some $\eta\in \Omega_F$.  There is a useful characterization of canonical divisors: $D\in \mathrm{Div}(F)$ is canonical if and only if 
\begin{equation}\label{charofcanondiv}
	\deg(D)=2g(F)-2\quad\text{and}\quad \dim_{\ff_q}\rr(D)\geq g(F),
\end{equation}
where $g(F)$ is the genus of $F$ and $\rr(D)$ is the Riemann-Roch vector space over $\ff_q$ defined as
\[\rr(D)=\{z\in F\setminus\{0\}:(z)+D\geq 0\}\cup\{0\}.\] 
Finally given an extension $F'/F$ of function fields and a Weil differential $\eta\in \Omega_F$ there is a Weil differential of $F'$, called the {\it cotrace} of $\eta$ in $F'/F$ and denoted as
\[\mathrm{Cotr}_{F'/F}(\eta),\]  
whose divisor satisfies
\[\bigl(\mathrm{Cotr}_{F'/F}(\eta)\bigr)=\bigl(\mathrm{Con}_{F'/F}(\eta)\bigr)+\mathrm{Diff}(F'/F).\]
Let $D\in \mathrm{Div}(F)$ and $D'\in\mathrm{Div}(F')$ be two divisors such that $D$ is canonical. We say that $D'$ is the {\it canonical lifting} of $D$ to $F'$ if $D'$ is not only a canonical divisor of $F'$ but also the lifting of $D$ to $F'$.

\subsection*{AG-codes and iso-dual AG-codes} \label{sec:isodual}
Let us recall now the definition of an AG-code. To construct an AG-code over $\fq$, a function field $F$ over $\fq$ and a couple of disjoint divisors of $F$ are needed. More precisely, for a function field $F$ over $\fq$,  we consider two disjoint divisors $D$ and $G$ of $F$ such that 
	$$D=P_1+\dots+P_n,$$ 
where $P_1,\ldots,P_n$ are pairwise distinct rational (i.e.\@ of degree 1) places of $F$. The AG-code $C_{\mathcal{L}}(D, G)$ over $F$ is the linear code over $\fq$ defined as
\begin{equation} \label{defAGcode}
	C_{\mathcal{L}}(D, G)=\{ (z(P_1),\ldots, z(P_n))\, :\, z\in \mathcal{L}(G)\} \subseteq \fq^n.
\end{equation}

Let $\cC \subseteq \fq^n$ be a linear code and ${\bf x} =(x_1, \ldots, x_n) \in (\fq^\ast)^n$. 
The set
	$$ {\bf x} \cdot \cC :=\{ (x_1 c_1, \ldots, x_n c_n) : (c_1, \ldots, c_n) \in \cC \} $$
is clearly another linear code over $\fq$. We notice that the codes $\cC$ and ${\bf x} \cdot \cC$ have the same length, dimension and minimum distance. 
This gives rise to the following notion of equivalence of linear codes as given in Definition 2.2.13 of \cite{Stichbook09}. 
\begin{definition}\label{defequiv}
	Two linear codes $\cC_1$ and $\cC_2$ over $\fq$ are \textit{equivalent} if $\cC_2={\bf x} \cdot\cC_1$ for some ${\bf x} \in (\fq^\ast)^n$. In this case we write $\cC_1 \sim_{\bf x} \cC_2$, or simply $\cC_1 \sim \cC_2$ when the distinction of the vector $\bf x$ is not necessary.
\end{definition}

The dual code $\cC^\perp$ of a linear  code $\cC$ over $\fq$  is the orthogonal complement of $\cC$ in $\fq^n$ with the standard inner product of $\fq^n$. It was shown in Corollary 2.6 of \cite{Stich1988} (see also Proposition 2.2.10 of \cite{Stichbook09}) that the dual of an AG-code $C_\rr(D,G)$ is again an AG-code $C_\rr(D,H)$ where $H=D-G+(\eta)$ for some $\eta\in \Omega_F$.

Recall that a linear code $\cC$ is said to be  {\em self-orthogonal} if $\cC \subset \cC^\perp$ and {\em self-dual} if $\cC=\cC^\perp$. The notion of equivalence of linear codes allows to define the following class of linear codes containing the class of self-dual codes.
\begin{definition}\label{defisodual}
	Let ${\bf x} \in (\fq^*)^n$. A code $\cC \subset \fq^n$ is called \textit{$\bf x$-iso-dual} if $\cC^\perp=\bf x\cdot \cC$. We will simply speak of \textit{iso-dual} codes when there is no need to specify the vector ${\bf x}$.
\end{definition}

Clearly the class of iso-dual codes contains the class of the self-dual codes, that is, $\bf x$-iso-dual codes are just self-dual codes with  ${\bf x}=(1,1,\ldots,1)$. 
We observe that the definition of self-dual code is a particular case of iso-dual codes where the equivalence of codes is given in terms of the monomial equivalence of codes, that is codes which are isometric to their duals via a permutation of coordinates and multiplying certain coordinates by non-zero constants (see, for instance, \cite{KL2017}).

\section{Lifting canonical divisors of the form $2G-D$ in a tower of function fields} \label{sec:lifting}
Let $F'/F$ be an extension of algebraic function fields over $\fq$. We recall a construction given recently in \cite{CPQT2024} (see also 
\cite{CPQT2024corr}): given two disjoint divisors $D$ and $G$ of $F$ such that $2G - D$ is equivalent 
to a canonical divisor of $F$ we were able to construct, under certain conditions on the ramification behavior in the extension $F'/F$, two disjoint divisors $\tilde{D}$ and $\tilde{G}$ of $F'$ such that $2\tilde{G}-\tilde{D}$ is a canonical divisor of $F'$ and also the divisor $2\tilde{G}-\tilde{D}$ of $F'$  is the lifting to $F'$ of  the canonical divisor  $2G-D$ of $F$.
 
 The conditions on $F'$ and $F$ allowing us to prove the above result are given in the next theorem, whose corrected proof can be found in 
\cite{CPQT2024corr} and the references therein.

\begin{theorem}[\cite{CPQT2024, CPQT2024corr}] \label{teolevantado}
    Let $F'/F$ be a finite separable extension of degree $m \geq 2$ of function 
    fields over $\mathbb{F}_q$. 
    Let $n, r \in \mathbb{N}$ with $n$ even, and suppose that 
    $A=\{P_1, \dots, P_n\}$ and $B=\{Q_1, \dots, Q_r\}$ are two disjoint sets of 
    places of $F$ such that:
    \begin{enumerate}[$(a)$]
        \item $P_1, \dots, P_n$ are rational places and $P_i$ splits completely 
        in $F'/F$ for $1 \leq i \leq n$,
        \item the extension $F'/F$ is unramified outside the set 
        $\{Q_1, \dots, Q_r\}$, and
        \item $d(R|Q_i)$ is even for each $1 \leq i \leq r$ and each place $R$ 
        of $F'$ lying over $Q_i$.
    \end{enumerate}
    Let $(\beta_1, \ldots, \beta_r) \in \mathbb{Z}^r$ be a non-zero $r$-tuple 
    and consider the disjoint divisors of $F$
    $$D = \sum_{i=1}^n P_i \qquad \text{and} \qquad G = \sum_{i=1}^r \beta_i Q_i,$$
    with $\deg(G) > 0$. Let $\tilde{D}$ and $\tilde{G}$ two divisors of $F'$ defined as
    \[\tilde{D} = \mathrm{Con}_{F'/F}(D) \qquad \text{and} \qquad 
    \tilde{G} = \mathrm{Con}_{F'/F}(G) + \tfrac{1}{2}\mathrm{Diff}(F'/F).\] 
    If $2G - D$ is a canonical divisor of $F$ then $2\tilde{G}-\tilde{D}$ is a canonical divisor of $F'$. More precisely if $2G-D=(\eta)$ for some Weil differential $\eta\in \Omega_F$ then $2\tilde{G}-\tilde{D}=(\tilde{\eta})$ where $\tilde{\eta}=\mathrm{Cotr}_{F'/F}(\eta)\in \Omega_{F'}$ and $2\tilde{G}-\tilde{D}$ is the canonical lifting of $2G-D$ to $F'$. In particular if we consider the AG-codes $C_{\mathcal{L}}(D, G)$ and $C_\rr(\tilde{D},\tilde{G})$ we have that both AG-codes are iso-dual. 
\end{theorem}
 
Under the hypothesis of Theorem \ref{teolevantado} we write
\[\mathfrak{L}_{F'/F}(C_\mathcal{L}(D,G))= C_{\mathcal{L}}(\tilde{D}, \tilde{G}),\]
 and with this notation we have that $\mathfrak{L}_{F'/F}(C_\mathcal{L}(D,G))$ is an AG-code over $F'$ which is called the \emph{lifting} of $C_{\mathcal{L}}(D,G)$ to $F'$.
\begin{remark}
	The condition of having $n\geq 2$ even in Theorem \ref{teolevantado} can not be removed. This is because if $g$ is the genus of $F$ and $2G-D$ is a canonical divisor of $F$ then, by  \eqref{charofcanondiv}, we have that 
	\[2g-2=\deg(2G-D)=2\deg(G)-n,\]
	and thus $n$ must be even. 
\end{remark}

The following estimates of the parameters of the lifting of an AG-code are straightforward to obtain, as can be seen in \cite{CPQT2024}.
\begin{corollary}[\cite{CPQT2024}] \label{parameters}
	Under the same conditions of Theorem \ref{teolevantado}, if $g_{_{F'}}$ is the genus of $F'$ and $nm> 2g_{_{F'}}-2$ 
	then $\mathfrak{L}_{F'/F}(C_\cL(D,G))$ is  an $[nm, \frac 12nm, d]$-code
	over $\fq$ where  
	$$ d \ge \tfrac 12nm -g_{_{F'}}+1. $$ 
\end{corollary}

We will use Theorem \ref{teolevantado} to lift a canonical divisor of the form $2G-D$ to each step of certain towers of function fields over a finite field. Let us recall first that a {\it tower} $\mathcal{T}$  of function fields  over a finite field $\fq$ is an infinite sequence of function fields $\mathcal{T}= \{F_i\}_{i \ge 0}$ such that every field extension $F_i/F_{i-1}$ is finite and separable, and such that $\fq$ is the full field of constants of every $F_i$. To avoid trivial situations, it is also required that  $g(F_i)\rightarrow\infty$ as $i \rightarrow \infty$. Towers of function fields have been found useful  for arithmetic purposes as well as  for some applications in Coding theory or Cryptography (see, for example, \cite{Stichbook09} and \cite{NX01}). For these applications two important sets of places of the first step $F_0$ of a tower $\mathcal{T}=\{F_i\}_{i \ge 0}$ are always important, namely the \textit{splitting locus} $\mathrm{Split}(\mathcal{T}/F_0)$ of $\mathcal{T}$ over $F_0$, which is the subset of $ \mathbb{P}_{F_0}$ defined as
\[ \mathrm{Split}(\mathcal{T}/F_0)=\{P\in \mathbb{P}_{F_0}: \text{$\deg P=1$ and $P$ splits completely in each $F_i$}\}, \]
and the \textit{ramification locus} $\mathrm{Ram}(\mathcal{T}/F_0)$ of $\mathcal{T}$ over $F_0$, which is the subset of $ \mathbb{P}_{F_0}$ defined as
\[ \mathrm{Ram}(\mathcal{T}/F_0)=\{P\in \mathbb{P}_{F_0}:\text{$P$ ramifies in  $F_i$ for some $i>0$}\}. \]
We show now that Theorem \ref{teolevantado} will yield a systematic way to propagate a canonical divisor of the form $2G-D$ along a tower of function fields over a finite field.
\begin{theorem}\label{liftingcandivtower}
   Let $\mathcal{T} = \{F_i\}_{i=0}^\infty$ be a tower of function 
    fields over $\fq$ and let $m_i=[F_{i+1}:F_i]$ for each $i\geq 0$. Let $n \geq 2$ be an even integer and suppose that $\{P_1,\ldots,P_n\}$ and $\{Q_1,\ldots,Q_r\}$ are two disjoint set places of $F_0$ such that
    \[\{P_1,\ldots,P_n\}\subset \mathrm{Split}(\mathcal{T}/F_0)\quad\text{and}\quad  \mathrm{Ram}(\mathcal{T}/F_0)\subset \{Q_1,\ldots,Q_r\}.\]
     Assume also  that for each index $i\geq 0$ all the different exponents defined by the places in the extension $F_{i+1}/F_i$ are even.

  Let us consider two disjoint divisors $D_0$ and $G_0$ of $F_0$ such that $2G_0-D_0$ is a canonical divisor. Then
  for each index $i\geq 1$ the divisors of $F_i$ defined 
    	as
    	\[ D_i = \mathrm{Con}_{F_i/F_0}(D_0) \quad \text{and} \quad 
    	G_i = \mathrm{Con}_{F_i/F_0}(G_0) + 
    	\tfrac{1}{2}\mathrm{Diff}(F_i/F_0),\]
    	are disjoint divisors of $F_i$ such that the divisor $2G_i - D_i$ is the canonical lifting to $F_i$ of the divisor $2G_0-D_0$. 
    
    \end{theorem}
\begin{proof}
    
    Let $i\geq 1$. We clearly have that $D_i$ and $G_i$ are disjoint divisors because the places in the support of the divisor $D_i$ lie over the places of the set $\{P_1,\ldots,P_n\}$ and the places in the support of the divisor $G_i$ lie over the places of the set $\{Q_1,\ldots,Q_r\}$ and $\{P_1,\ldots,P_n\}\cap \{Q_1,\ldots,Q_r\}=\emptyset$. To finish the proof it suffices to check that the extension $F_i/F_0$ satisfies the conditions $(a)$, $(b)$ and $(c)$ of Theorem \ref{teolevantado} with 
    \[A=\{P_1,\ldots,P_n\}\quad\text{and}\quad B=\{Q_1,\ldots,Q_r\}.\]
    Clearly conditions $(a)$ and $(b)$ of Theorem \ref{teolevantado} hold for the extension $F_i/F_0$ by the ramification properties of the places of $\mathrm{Split}(\mathcal{T}/F_0)$  and $\mathrm{Ram}(\mathcal{T}/F_0)$. We proceed by induction to prove that  condition $(c)$ of \ref{teolevantado} holds for the extension $F_i/F_0$ for any index $i\geq 1$: the case $i=1$ is just an hypothesis in our theorem.  Suppose now that for any place $S$ of $F_i$ lying over a place $Q\in \mathrm{Ram}(\mathcal{T}/F_0)$ the different exponent $d(S|Q)$ is even. Let $R$ be a place of $F_{i+1}$ lying over a place $Q\in \mathrm{Ram}(\mathcal{T}/F_0)$. If $S=R\cap F_i$ then $S|Q$ and by \cite[Collorary 3.4.12]{Stichbook09} we have that
    \[d(R|Q)=e(R|S)d(S|Q)+d(R|S).\]
    Since $d(S|Q)$ is even by inductive hypothesis and $d(R|S)$ is even by hypothesis, we see that $d(R|Q)$ is also even. Therefore the extension $F_i/F_0$ also satisfies the condition $(c)$ of Theorem~\ref{teolevantado}, with the set of places $A$ and $B$ as above, for any index $i\geq 1$. Thus item $(iii)$ holds and we are done. 
\end{proof}
A direct coding theoretic consequence of Theorem \ref{liftingcandivtower} is the following
\begin{corollary}\label{liftisodualtower}
	Under the conditions of Theorem \ref{liftingcandivtower} we have  that the lifting $\mathfrak{L}_{F_i/F_0}(C_\mathcal{L}(D_0,G_0))$ of $C_\mathcal{L}(D_0,G_0)$ to $F_i$ is  an iso-dual code for any index $i\geq 0$. 
\end{corollary}

In order to be able to use Theorem \ref{liftingcandivtower} we have to start with two divisors $D_0$ and $G_0$ of $F_0$ such that $2G_0-D_0$ is a canonical divisor. This is  a problem that can be very difficult to solve in general  but, fortunately, when $F_0$ is a rational function field canonical divisors are easily identified by degree considerations and a simple condition ensures that $2G_0-D_0$ is a canonical divisor. The following result covers the most common setting in applications, namely when the first step in a tower  is a rational function field over a finite field. 
\begin{corollary}\label{liftingfromratfuncfield}
Under the assumptions of Theorem~\ref{liftingcandivtower}, suppose that
$F_0$ is the rational function field over $\ff_q$. If $D_0$ and $G_0$ are disjoint divisors of $F_0$ such that
\[
\deg(G_0)=\frac{n-2}{2},
\]
then for each index $i\geq 1$ the divisors of $F_i$ defined 
	as
	\[ D_i = \mathrm{Con}_{F_i/F_0}(D_0) \quad \text{and} \quad 
	G_i = \mathrm{Con}_{F_i/F_0}(G_0) + 
	\tfrac{1}{2}\mathrm{Diff}(F_i/F_0),\]
	are disjoint divisors of $F_i$ such that the divisor $2G_i - D_i$ is the canonical lifting to $F_i$ of the divisor $2G_0-D_0$.  
\end{corollary}
\begin{proof}
	In view of Theorem \ref{liftingcandivtower} we see that all we need to show is that $2G_0-D_0$ is a canonical divisor of the rational function field $F_0$. This follows directly from \cite[Proposition 1.6.2]{Stichbook09} because $g(F_0)=0$ and
	\[\deg(2G_0-D_0)=2\deg(G_0)-\deg(D_0)=n-2-n=-2.\]
\end{proof}

\section{Good sequences of iso-dual AG-codes from towers} \label{sec: goodsequences}

A natural question is when the sequences of iso-dual AG-codes given in Corollary \ref{liftisodualtower} are good sequences in a precise sense defined below.  We will see in this section that Corollary \ref{liftisodualtower} can be used in some good recursive towers of function fields and, as stated in the Introduction, in this way we will prove the existence of good sequences of iso-dual AG-codes in a simpler way than the construction of good self-dual AG-codes given in \cite{Sch89}, \cite{Stich06} and \cite{BS2019} with the advantage of using recursive towers function fields which is the most desirable situation.

First we recall what  a good sequence of linear codes is. Let $\cC$ be a linear $[n, k, d]$-code over $\fq$. The ratios 
$$R(\cC ) = \tfrac{k(\cC )}{n(\cC)}=\tfrac{k}{n} \qquad \text{ and } \qquad \delta (\cC ) = \tfrac{d(\cC)}{n(\cC)}=\tfrac{d}{n},$$ 
denote the \textit{information rate}
and the \textit{relative minimum distance} of the code, respectively. An infinite sequence $\{\mathcal{C}_i\}_{i\ge 1}$ of $[n_i,k_i,d_i]$-codes over $\ff_q$ is called \textit{good} (or {\it asymptotically good}) if the information rate $R(\mathcal{C}_i) = k_i/n_i$ and the relative minimum distance $\delta(\mathcal{C}_i) = d_i/n_i$ of $\mathcal{C}_i$ satisfy 
\[ \liminf_{i \rightarrow \infty} R(\mathcal{C}_i) > 0 \qquad \text{and} \qquad \liminf_{i \rightarrow \infty} \delta(\mathcal{C}_i)>0. \] 

An infinite sequence $\{\mathcal{C}_i\}_{i \ge 1}$ of $[n_i,k_i,d_i]$-codes over $\ff_{\ell^2}$ is said to {\it attain} the
TVZ-bound if, for  $0< \delta\leq 1-(\ell-1)^{-1}$
, we have
$$ \liminf_{i \rightarrow \infty} R(\mathcal{C}_i) \geq 1- \frac{1}{\ell-1}-\delta \qquad \text{and} \qquad \liminf_{i \rightarrow \infty} \delta(\mathcal{C}_i) \geq  \delta. $$ 

It turns out that towers of function fields over a finite field $\fq$ are suitable to construct asymptotically good sequences of linear codes.  We refer to Chapter 7 of \cite{Stichbook09} for more details on definitions and proofs of the results mentioned below.

 The suitability of a tower $\mathcal{T}=\{F_i\}_{i\geq 0}$  over $\fq$ for these applications in Coding theory is determined by the following two  quantities associated to $\mathcal{T}$. The first quantity is the \textit{genus} $\gamma(\mathcal{T})$ of the tower $\mathcal{T}$ over $F_0$, which is defined as 
\begin{equation} \label{eq: genusT}
	\gamma(\mathcal{T}):=\lim_{i \rightarrow \infty}\frac{g(F_i)}{[F_i:F_0]}.
\end{equation} 
The second quantity is the 
\textit{splitting rate} $\nu(\mathcal{T})$ of the tower $\mathcal{T}$ over $F_0$, which is defined as
\[ \nu(\mathcal{T}):= \lim_{i \rightarrow \infty}\frac{N(F_i)}{[F_i:F_0]},\]
where $N(F_i)$ denotes the number of rational places (or places of degree one) of $F_i$.

These two quantities determine the asymptotic behavior of a tower of function fields over $\ff_q$ as follows. The \textit{limit} of a tower $\mathcal{T}$ over $\ff_q$ is defined as
\[ \lambda(\mathcal{T}):= \lim_{i \rightarrow \infty} \frac{N(F_i)}{g(F_i)}. \]
It can be proved that $0<\gamma(\mathcal{T})\leq \infty$ and $0\leq \nu(\mathcal{T})<\infty$. Thus, since we have the equation $$\lambda(\mathcal{T})=\frac{\nu(\mathcal{T})}{\gamma(\mathcal{T})},$$ it turns out that $0\leq \lambda(\mathcal{T})<\infty$ for any tower $\mathcal{T}$ over $\ff_q$ and the good or bad asymptotic behavior of $\mathcal{T}$ is decided according to whether or not its limit is zero. 
More precisely, a tower $\mathcal{T}$ of function fields over  $\ff_q$ is said to be \textit{(asymptotically) good} if $\lambda(\mathcal{T})>0$. Otherwise, that is $\lambda(\mathcal{T})=0$, the tower is said to be \textit{(asymptotically) bad}. 

Roughly speaking, a tower $\mathcal{T}=\{F_i\}_{i \ge 0}$ over $\ff_q$ is good if $N(F_i)$ grows faster than  the genus of $F_i$ as $i$ goes to infinity. These are the kind of towers of interest in the above mentioned applications. It is well known that if $\mathcal{T}$ is a tower over $\ff_q$, then
\begin{equation}\label{ihara}
	\lambda(\mathcal{T})\leq\sqrt{q}-1.
\end{equation}
A tower $\mathcal{T}$ of function fields  over $\ff_q$ is called {\it optimal} if equality in \eqref{ihara} occurs.

The genus $\gamma(\mathcal{T})$ and the splitting rate $\nu(\mathcal{T})$ of a tower $\mathcal{T}=\{F_i\}_{i \ge 0}$ of function fields over $\fq$ are usually studied in terms of
the splitting locus $\mathrm{Split}(\mathcal{T}/F_0)$ and the ramification locus $\mathrm{Ram}(\mathcal{T}/F_0)$  of $\mathcal{T}$ over $F_0$ defined in the previous section. In many cases  the finiteness of $\mathrm{Ram}(\mathcal{T}/F_0)$ together with the non-emptyness of $\mathrm{Split}(\mathcal{T}/F_0)$ are the two conditions allowing to prove that the tower $\mathcal{T}$ is good.

All the examples we will present next are given in terms of the so called \textit{recursive towers} of function fields, which are the most useful in the applications. Recursive means that the first function field in the tower is a rational function field, that is  $F_0=\fq(x_0)$ for some transcendental element $x_0$ over $\fq$, and all the others function fields in the sequence can be defined recursively for every $i\geq 0$ as $F_{i+1}=F_{i}(x_{i+1})$, where $\{x_i\}_{i=0}^{\infty}$ is a sequence of transcendental elements over $\fq$ satisfying an equation of the form $f(x_{i},x_{i+1})=0$, for some appropriate bivariate polynomial $f \in \fq[x,y]$.

From Corollary \ref{liftisodualtower} we can give  conditions to obtain asymptotically good sequences of iso-dual AG-codes over $\ff_q$ from towers of function fields.

\begin{theorem}\label{isontowers}
Let $n\geq 2$ be an even integer and let $\mathcal{T}=\{F_i\}_{i \ge 0}$ be a tower of function fields over $\fq$ of genus 
	\[ \gamma(\mathcal{T}) \le n (\tfrac 12 -\delta), \] 
for some $0<\delta< \frac 12$. Let $r$ be a positive integer and suppose that   $\{P_1,\ldots,P_n\}$ and  $\{Q_1,\ldots,Q_r\}$ are two disjoint sets of places of $F_0$ such that
\begin{enumerate}[$(a)$]
	\item $P_1,\ldots, P_n$ are rational and $\{P_1,\ldots,P_n\}\subset \mathrm{Split}(\mathcal{T}/F_0)$.\sk
	
	\item 	$\mathrm{Ram}(\mathcal{T}/F_0)\subset \{Q_1,\ldots,Q_r\}$.
	  \sk
	
	\item For each index $i\geq 1$ the different exponents defined by the places of the extension $F_{i+1}/F_i$ are even.
\end{enumerate}
Consider the disjoint divisors of $F_0$ 
	\[ D_0=P_1 + \cdots + P_n \qquad \text{and} \qquad G_0=\sum_{1\le j \le r} \b_jQ_j, \] 
for some integers $\b_1,\ldots, \b_r$, and suppose that the divisor $2G_0 - D_0$ is a canonical divisor of $F_0$. Then the sequence of liftings $\{\mathcal{C}_i\}_{i\geq 0}$, where			
	$$ \mathcal{C}_i= \mathfrak{L}_{F_i/F_0}(C_\cL(D_0,G_0)), $$
is an asymptotically good sequence of iso-dual AG-codes over $\fq$. More precisely, for each index $i\geq 0$ the linear code
$\mathcal{C}_i$ 
is an 
iso-dual AG-code over $\fq$ of length $n_i=n[F_i:F_0]$  such that 
\begin{equation}\label{isorates}
\liminf_{i\rightarrow\infty} R(\mathcal{C}_i)= \tfrac{1}{2} \qquad \text{and} \qquad \liminf_{i\rightarrow\infty} \delta(\mathcal{C}_i)\geq \delta>0,
\end{equation}
where $k_i$ and $d_i$ are, respectively, the dimension and the minimum distance of $\mathcal{C}_i$ 
for each $i> 0$.  
\end{theorem}

\begin{proof}
We are in the conditions of Theorem \ref{liftingcandivtower} so that from Corollary \ref{liftisodualtower} we have that for each index $i\geq 0$ the lifting $\mathcal{C}_i=\mathfrak{L}_{F_i/F_0}(C_\cL(D_0,G_0))$ of $C_\rr(D_0, G_0)$ to $F_i$ is an iso-dual AG-code of length $n_i=n[F_i:F]$. 

Now let $g_i$ be the genus of $F_i$ and let $m_i=[F_i:F_0]$ for each $i\geq 0$.  The sequence
	\[ \Big\{ \frac{g_i-1}{m_i} \Big \}_{i\geq 0} \]
is monotonically increasing (see, for example, Lemma 7.2.3 of \cite{Stichbook09}) and
	\[ \lim_{i\rightarrow\infty}\frac{g_i-1}{m_i}=\lim_{i\rightarrow\infty} \frac{g_i}{m_i} = \gamma(\mathcal{T}/F_0) \le n (\tfrac{1}{2}-\delta) < \tfrac 12 n. \]
Then we have that there exists $i_0 >0$ such that
	\[ \frac{g_i-1}{m_i} \le n (\tfrac{1}{2}-\delta) < \tfrac{1}{2} n, \]
for each $i\geq i_0$, and thus $n_i=nm_i>2g_i-2$ for $i\geq i_0$. From Corollary \ref{parameters} we see that
	\[ k_i = \tfrac 12 n_i \qquad \text{and} \qquad d_i \ge \tfrac 12 (n_i-2g_i+2), \]
so that
\[ \liminf_{i\rightarrow\infty}R(\mathcal{C}_i) = \frac{1}{2}\quad\text{and}\quad\liminf_{i\rightarrow\infty}\delta(\mathcal{C}_i) \ge \liminf_{i\rightarrow\infty} \Big(\frac{1}{2}-\frac{g_i-1}{n_i}\Big) \ge \delta>0,\]
as claimed.
\end{proof}

\begin{remark}
		The condition $\gamma(\mathcal{T})\leq n(\tfrac1 2-\delta)$ in Theorem \ref{isontowers} implies that $\mathcal{T}$ has to be an asymptotically good tower over $\ff_q$ with limit $\lambda(\mathcal{T})> 2$. This is so because (with the notation used in Theorem \ref{isontowers})
		$N(F_i)\geq n[F_i:F_0]$ so that
		\[\lambda(\mathcal{T})=\lim_{i \rightarrow \infty}\frac{N(F_i)}{g(F_i)}\geq n\lim_{i \rightarrow \infty}\frac{[F_i:F_0]}{g(F_i)}=\frac{n}{\gamma(\mathcal{T})}\geq \left(\tfrac{1}{2}-\delta\right)^{-1}> 2.\]
		
		In particular from the upper bound \eqref{ihara} for $\lambda(\mathcal{T})$ we have that  
		\[2<\lambda(\mathcal{T})\leq \sqrt{q}-1,\]
		so that $q>9$ and thus	Theorem \ref{isontowers} holds for finite fields $\ff_q$ with $q\geq 11$.
\end{remark}

In order to use Theorem \ref{isontowers} in a tower $\mathcal{T}=\{F_i\}_{i \ge 0}$ to obtain good sequences of iso-dual AG codes, we need to start with divisors $D_0$ and $G_0$ of $F_0$ such that $2G_0-D_0$ is a canonical divisor. As we have shown in Corollary \ref{liftingfromratfuncfield} there is a simple condition to guarantee this when the first step $F_0$ of the tower $\mathcal{T}$ is a rational function field. More precisely we have the following  

\begin{corollary}\label{isoratexamples}
Suppose that $\mathcal{T}=\{F_i\}_{i \ge 0}$ is a tower of function fields over $\ff_q$ satisfying the conditions  
of Theorem \ref{isontowers}, where $F_0=\ff_q(x)$. Consider the disjoint divisors of $F_0$
	\[ D_0=P_1 + \cdots + P_n \qquad \text{and} \qquad G_0=\sum_{1 \le j \le r} \b_jQ_j \] 
and the associated AG-code $C_{\cL}(D_0,G_0)$. 
If the integers $\b_1,\ldots,\b_r$ are such that
	\[ \sum_{1\le i \le r} \b_i\deg Q_i = \tfrac 12 (n-2),\]
then the sequence of liftings $\{ \mathfrak{L}_{F_i/F_0}(C_\cL(D_0,G_0))\}_{i\geq 0}$ defines an asymptotically good sequence of iso-dual AG-codes over $\fq$ with rates satisfying \eqref{isorates}.
\end{corollary}

\begin{proof}
From Corollary \ref{liftingfromratfuncfield} we have that the lifting $\mathcal{C}_i= \mathfrak{L}_{F_i/F_0}(C_\cL(D_0,G_0))$ of $C_\rr(D_0,G_0)$ to $F_i$ is an iso-dual AG-code over $\ff_q$. The conclusion follows at once from Theorem \ref{isontowers}.
\end{proof}

In what follows, we present two explicit constructions of sequences of asymptotically good iso-dual codes from well known  recursive towers of function fields.

\subsection{A sequence of iso-dual AG-codes attaining the TVZ-bound from a recursive tower} 
\label{goodisodualq2} 

It is well known that there exists a family of self-dual AG-codes attaining the TVZ bound. In characteristic two the existence of self-orthogonal and self-dual codes attaining the TVZ bound was first proved in \cite{M2002}. In \cite{Stich06} the same conclusion is achieved in any characteristic by using the Galois closure of an asymptotically optimal tower of function fields over $\ff_{q^2}$ presented in  \cite{GS96}. Since the information rate $R$ in the case of self-dual and iso-dual codes is 
	$$R= \tfrac 12,$$ 
having a sequence $\{\mathcal{C}_i\}_{i\ge 0}$ of $[n_i,k_i,d_i]$-self-dual codes (or, more generally, iso-dual codes) over $\ff_{q^2}$ attaining the TVZ bound means that 
	$$\delta = \tfrac 12 - (q-1)^{-1}$$ 
and so the relative minimum distance $\delta(\mathcal{C}_i)$ of $\mathcal{C}_i$ must satisfy
\begin{equation}\label{mindisisodual}
	\liminf_{i \rightarrow \infty} \delta(\mathcal{C}_i)\geq \delta = \frac 12 - \frac{1}{q-1}.
\end{equation}
   
 We are going to show that we can use Corollary \ref{isoratexamples} with the asymptotically optimal 
 recursive tower of function fields over $\mathbb{F}_{q^2}$ given in \cite{GS96} to get, in a simpler way, a sequence of iso-dual AG-codes attaining the TVZ bound. 

First, we recall some results about the recursive tower given in \cite{GS96} needed to use Corollary~\ref{isoratexamples}. Let $q$ be prime power and let $\mathcal{T}=\{F_i\}_{i \ge 0}$ be  the tower of function fields over $\mathbb{F}_{q^2}$ recursively defined by  $F_0=\mathbb{F}_{q^2}(x_0)$ and  $F_i=F_{i-1}(x_i)$, where for each $i\geq 1$ we have
\begin{equation}\label{eqtower}
x_i^q+x_i=\frac{x^q_{i-1}}{1+x^{q-1}_{i-1}}.
\end{equation}

For $\alpha \in \fqs$, we denote by $P_\alpha$ the only zero of $x_0-\alpha$ in $F_0$, and by $P_\infty$ the only pole of $x_0$ in $F_0$. Consider now the set 
$$\Omega = \{\alpha\in \mathbb{F}_{q^2}:\alpha^q+\alpha=0\}.$$ 

It was shown in Lemmas 3.3, 3.5 and 3.9 of \cite{GS96} that  
\begin{align*}
	& \textrm{Split}(\mathcal{T}/F_0)=\{P_\alpha\in \mathbb{P}_{F_0}:\alpha\notin \Omega\}, \\[1mm] 
	& \textrm{Ram}(\mathcal{T}/F_0)\subset\{P_\alpha\in  \mathbb{P}_{F_0}:\alpha\in \Omega\}\cup \{P_\infty\}.
\end{align*}	
Also, from Remark 3.8 of \cite{GS96}, the genus of $F_n$ is given by
	$$ g(F_n)=\left\{\begin{array}{ll}
	(q^{\frac{n+1}2}-1)^2, 				 & \quad \text{if $n$ is odd}, \\[1.5mm]
	(q^{\frac{n+2}2}-1)(q^{\frac n2}-1), & \quad \text{if $n$ is even}.
\end{array}\right. $$ 
Since $[F_n:F_0]=q^n$, from \eqref{eq: genusT} we conclude that 
	$$\gamma(\mathcal{T}) \le q.$$
Moreover, the cardinality of the splitting set is 
	$$|\textrm{Split}(\mathcal{T}/F_0)|=q^2-|\Omega|=q^2-q.$$ 
\begin{proposition} \label{prop: TVZ}
There is a sequence $\{\mathcal{C}_i\}_{i \ge 1}$ of iso-dual AG-codes over $\ff_{q^2}$ attaining the TVZ bound constructed from a recursive tower of function fields over $\ff_{q^2}$ with $q\geq 4$. 
\end{proposition}

\begin{proof}
Consider the tower $\mathcal{T}=\{F_i\}_{i \ge 0}$  of function fields over $\mathbb{F}_{q^2}$ recursively defined by the equation \eqref{eqtower}. 

For $n=q^2-q$ and $ \delta=\tfrac{1}{2}-\tfrac{1}{q-1}$, 
we have not only that $n\geq 2$ is even but also that
	\[ \gamma(\mathcal{T}) \le q \le n \left(\tfrac{1}{2}-\delta\right), \] 
for any  $q\geq 4$. Now, for any place $P$ of $F_n$ lying over $P_\alpha$ with $\alpha\in \Omega$ or over $P_\infty$, we have from Lemmas~3.3 and 3.5 of \cite{GS96} that either $P$ is unramified in $F_n/F_{n-1}$ or else $P$ is totally ramified in $F_n/F_{n-1}$ with different exponent $d(P|Q)=2(q-1)$, where $Q=P\cap F_{n-1}$.   By using the well known formula (see Corollary 3.4.12 of \cite{Stichbook09}) 
\[d(P|R)=e(P|Q)d(Q|R)+d(P|Q),\]
where $R=Q\cap F_{n-2}$, we  conclude that both different exponents $d(P|P_\infty)$ and $d(P|P_\alpha)$ are even for any $\alpha\in \Omega$. 

Consider the following divisors in the rational function field $F_0$,
	\[D_0=\sum_{\alpha\notin \Omega}P_\alpha \qquad \text{and} \qquad 
		G_0=\sum_{\alpha\in\{\infty\}\cup \Omega}\beta_\alpha P_\alpha,\]
where the integers $\beta_\alpha$ are such that 
	\[ \sum_{\alpha\in\{\infty\}\cup \Omega}\beta_\alpha = \tfrac12 (q^2-q-2). \]
Thus, since $q \geq 4$, we have from Corollary \ref{isoratexamples} that the sequence of liftings 
	$$\mathcal{C}_i = \mathfrak{L}_{F_i/F_0}(C_\cL(D_0,G_0))$$ 
is a good sequence of iso-dual AG-codes over $\mathbb{F}_{q^2}$ for $q\geq 4$ with the following rates given by \eqref{isorates}
	\[\liminf_{i\rightarrow\infty} R(\mathcal{C}_i) = \frac{1}{2} \qquad \text{and} \qquad \liminf_{i\rightarrow\infty}\delta(\mathcal{C}_i)\geq \frac{1}{2}-\frac{1}{q-1} .\] 
In particular we see that \eqref{mindisisodual} holds so that the sequence of iso-dual AG-codes $\{\mathfrak{L}_{F_i/F_0}(C_\cL(D,G))\}_{i>0}$ over $\ff_{q^2}$  attains the TVZ  bound for any prime power $q\geq 4$.	
\end{proof}

\subsection{An asymptotically good sequence of iso-dual codes in characteristic 2} 

In this section we make use the recursive tower of function fields $\mathcal{T}=\{F_i\}_{i\ge 0}$  given in 2005 by Bezerra, García and Stichtenoth \cite{BGS05} to produce good sequences of iso-dual codes over $\ff_{q^3}$ of even characteristic. We refer to this tower as the BGS-tower.

\begin{proposition} \label{prop: goodisodualF8}
For any integer $r\geq 2$ there exists a good sequence $\{\mathcal{C}_i\}_{i\ge 0}$ of iso-dual AG-codes  over $\ff_{8^r}$ constructed from a recursive tower of function fields over $\ff_{8^r}$ with rates  
\begin{equation}\label{eq: ratesF8}
	\liminf_{i\rightarrow\infty}R(\mathcal{C}_i)= \frac{1}{2} \qquad \text{and} \qquad 
	\liminf_{i\rightarrow\infty}\delta(\mathcal{C}_i) \ge \frac{1}{2}-\frac{2^r+2}{2(4^r-1)} .
\end{equation}   
\end{proposition}
\begin{proof}
We first recall the definition of the BGS-tower and some of their properties.	
Let $s\in \N$ and let $q$ be a prime power. Let $\mathcal{T}=\{F_i\}_{i \ge 0}$ be the 
tower over $\ff_{q^s}$ recursively defined by the equation
	\[ \frac{1-y}{y^q}=\frac{x^q+x+1}{x}. \]
That is, $F_0=\ff_{q^s}(x_0)$ and $F_i=F_{i-1}(x_i)$ where, for $i\ge 1$, we have 
	\[ \frac{1-x_i}{x_i^q}=\frac{x_{i-1}^q+x_{i-1}+1}{x_{i-1}}. \] 

It was proved in Proposition 1.3 of \cite{BGS05}  that $\mathrm{Ram}(\mathcal{T}/F_0)$ is a finite set and that  $P_\infty\in \mathrm{Ram}(\mathcal{T}/F_0)$, where $P_\infty$ is the only pole of $x_0$ in $F_0$. From  Theorem 2.9 of \cite{BGS05} we easily deduce that the genus of this tower over $\ff_{q^s}$ is
	\[ \gamma(\mathcal{T}) = \frac{q(q+2)}{2(q-1)}. \]
It was also proved in Theorem 3.2 of \cite{BGS05}  that when $s=3$ there are $n=q(q+1)$ rational places $P_{\alpha_1},\ldots,P_{\alpha_n}$ of $F_0$ splitting completely in the tower $\mathcal{T}$. Moreover, it was  shown in Propositions 2.5, 2.6, 2.7 and 2.8 of \cite{BGS05} that the different exponents are always even only when $q$ is even.

Thus, we can apply our method given in Theorem \ref{isontowers} by considering this tower $\mathcal{T}$ over cubic finite fields of even characteristic. For this tower $\mathcal{T}$ over $\ff_{q^3}$ with $q\geq 4$ even (that is $q=2^r$ with $r\geq 2$) by taking 
	$$\delta=\frac{1}{2}-\frac{q+2}{2(q^2-1)}, $$ 
we have not only that $0<\delta< \frac 12$, but also that
	$$\gamma(\mathcal{T}) = \frac{q(q+2)}{2(q-1)}= n(\tfrac{1}{2}-\delta).$$
Therefore by considering the following divisors 
	\[ D=\sum_{1\le i \le n} P_{\alpha_i} \qquad \text{and} \qquad G= \tfrac 12 (n-2) P_\infty \]
of $F_0$, we have 
from Corollary \ref{isoratexamples} that the sequence of liftings  $\{\mathfrak{L}_{F_i/F_0}(C_\cL(D,G))\}_{i>0}$ defines a good sequence of iso-dual AG-codes over $\mathbb{F}_{q^3}$ with rates given by \eqref{isorates}, hence \eqref{eq: ratesF8} holds.
\end{proof}

\section{Final remarks: successive liftings}\label{succliftings}
Suppose that we have a tower $\mathcal{T}=\{F_i\}_{i\ge 0}$ of function fields and a canonical divisor of $F_0$ of the form $2G_0-D_0$. Under the conditions of Theorem \ref{liftingcandivtower} we proved that in each step $F_i$ of the tower $\mathcal{T}$ we have a divisor of the form $2G_i-D_i$ which is the canonical lifting of $2G_0-D_0$  to $F_i$. A possible alternative approach could be to iterate Theorem \ref{teolevantado} step by step in the tower, that is, we first consider the canonical lifting $2G'_1-D'_1$ of $2G_0-D_0$ to $F_1$, then the canonical lifting $2G^{''}_2-D^{''}_2$ of $2G'_1-D'_1$ to $F_2$ and so on until we get a canonical divisor of $F_i$ of the form $2G-D$ where $D$ and $G$ are disjoint divisors of $F_i$. A natural question here is whether or not we arrive to the same canonical divisor of $F_i$. The answer is affirmative and to see this it suffices to check it in a finite tower of function fields of three steps. 
\begin{figure}
\begin{center}
	\begin{tikzpicture}
		\draw (-7,0)--(-7,3);
		
		\filldraw [white] (-7,0) circle (8pt);
		\node at (-7,0) {$F_0$};
		\filldraw [white] (-7,1.5) circle (8pt);
		\node at (-7,1.44) {$F_1$};
		\filldraw [white] (-7,3) circle (8pt);
		\node at (-7,3) {$F_2$};
		
		\draw (-5.5,0)--(-5.5,3);
		\filldraw [white] (-5.5,0) circle (8pt);
	\node at (-4.4,0) {{\small $(\eta_{_0})=2G_0-D_0$}};	
		\filldraw [white] (-5.5,3) circle (8pt);
		\node at (-4.6,3) {{\small $(\tilde{\eta})=2\tilde{G}-\tilde{D}$}};
		\node at (-3,1.5){{\small $\tilde{\eta}=\mathrm{Cotr}_{F_2/F_0}(\eta_{_0})$}};
		
		\draw (0,0)--(0,3);
		\filldraw [white] (0,0) circle (8pt);
		\node at (1.1,0) {{\small $(\eta_{_0})=2G_0-D_0$}};
		\filldraw [white] (0,1.5) circle (8pt);
		\node at (1.1,1.5) {{\small $(\eta_{_1})=2G_1-D_1$}};
		\filldraw [white] (0,3) circle (8pt);
		\node at (1.1,3) {{\small $(\eta_{_2})=2G_2-D_2$}};
		\node at (5,1.5){{\small $\eta_{_1}=\mathrm{Cotr}_{F_1/F_0}(\eta_{_0})$}};
		\node at (5,3){{\small $\eta_{_2}=\mathrm{Cotr}_{F_2/F_1}(\eta_{_1})$}};
	\end{tikzpicture}
	\caption{Succesive liftings vs one-step lifting}\label{succlifting}
\end{center}	
\end{figure}

More precisely suppose that in the finite tower (see Figure \ref{succlifting} below)  
	$$ F_0\subset F_1\subset F_2 $$ 
of function fields over $\ff_q$ the extension $F_2/F_0$ satisfies the conditions of Theorem \ref{teolevantado} with the rational places $P_1\ldots,P_n$ of $F_0$ splitting completely in $F_2$, the extension $F_2/F_0$ being unramified outside the places $Q_1,\ldots,Q_r$ and the different exponent $d(R|P)$ being even for any place $R$ of $F_2$ lying over a place $P$ of $F_0$. If there is a Weil differential $\eta_0\in \Omega_{F_0}$ such that $(\eta_{0})=2G_0-D_0$ then from Theorem \ref{teolevantado} we know that if $\tilde{\eta}=\mathrm{Cotr}_{F_2/F_0}(\eta_{_0})$ then
$\bigl(\tilde{\eta}\bigr)=2\tilde{G}-\tilde{D}$ where $\tilde{D}=\mathrm{Con}_{F_2/F_0}(D_0)$ and $\tilde{G}=\mathrm{Con}_{F_2/F_0}(D_0)+\tfrac{1}{2}\mathrm{Diff}(F_2/F_0)$. In particular the places $P_0,\ldots,P_n$ split completely in $F_1$ and the extension $F_1/F_0$ is unramified outside the places $Q_1,\ldots,Q_r$. Suppose that the different exponents $d(R|Q)$ and $d(Q|P)$ are even for any place $R$ of $F_2$ where $Q=R\cap F_1$ and $P=Q\cap F_0$. Then we can apply Theorem \ref{teolevantado} to the extensions $F_2/F_1$ and $F_1/F_0$ (see Figure \ref{succlifting} above): first we have that the Weil differential $\eta_{_1}=\mathrm{Cotr}_{F_1/F_0}(\eta_{_0})\in \Omega_{F_1}$ and the divisors of $F_1$
\[D_1=\mathrm{Con}_{F_1/F_0}(D_0)\quad\text{and}\quad G_1=\mathrm{Con}_{F_1/F_0}(G_0)+\tfrac{1}{2}\mathrm{Diff}(F_1/F_0),\] 
are such that $(\eta_{_1})=2G_1-D_1$. Now if $\{S_1,\ldots,S_m\}=\sop(D_1)$, then these places are rational places of $F_1$, they  form the overset of $P_1\ldots,P_n$ in $F_1$ and they split completely in $F_2$. On the other hand if $\{T_1,\ldots,T_k\}=\sop(G_1)$ then the extension $F_2/F_1$ is unramified outside the places $T_1,\ldots,T_k$. Then the extension $F_2/F_1$ satisfies the conditions of Theorem \ref{teolevantado} with the disjoint set of places $\{S_1,\ldots,S_m\}$ and $\{T_1,\ldots,T_k\}$ and since $2G_1-D_1$ is a canonical divisor of $F_1$, we have from Theorem \ref{teolevantado} that the Weil differential  $\eta_{_2}=\mathrm{Cotr}_{F_2/F_1}(\eta_{_1})\in \Omega_{F_2}$ and the divisors of $F_2$
\[D_2=\mathrm{Con}_{F_2/F_1}(D_1)\quad\text{and}\quad G_2=\mathrm{Con}_{F_2/F_1}(G_1)+\tfrac{1}{2}\mathrm{Diff}(F_2/F_1),\] 
are such that $(\eta_{_2})=2G_2-D_2$. We show now that
\[\eta_{_2}=\tilde{\eta}, \quad G_2=\tilde{G}\quad\text{and}\quad D_2=\tilde{D}.\]
Since the places $P_1,\ldots,P_n$ are rational places of $F_0$ which split completely in $F_2$, it is easy to see that
\[D_2=\mathrm{Con}_{F_2/F_1}(D_1)=\mathrm{Con}_{F_2/F_1}(\mathrm{Con}_{F_1/F_0}(D_0))=\mathrm{Con}_{F_2/F_0}(D_0)=\tilde{D}.\]
On the other hand by \cite[Proposition 3.4.11]{Stichbook09}
\[
\eta_2 = \mathrm{Cotr}_{F_2/F_1}(\eta_{_1})= \mathrm{Cotr}_{F_2/F_1}(\mathrm{Cotr}_{F_1/F_0}(\eta_{_0}))=\mathrm{Cotr}_{F_2/F_0}(\eta_{_0})=\tilde{\eta},
\]
so that
\[2G_2-\tilde{D}=2G_2-D_2=(\eta_{_2})=(\tilde{\eta})=2\tilde{G}-\tilde{D}.\]
Therefore $G_2=\tilde{G}$ and we are done.

In particular we see that lifting iso-dual AG codes step by step in a tower of function fields produces the same sequence of iso-dual AG-codes constructed in Theorem~\ref{isontowers} by lifting directly an iso-dual AG-code over the bottom field to each function field in the tower.

\section{Acknowledgements}
 
 The first author was partially supported by Despegue científico (PEDECIBA). The second and fourth authors were partially supported by CONICET, FONCyT, SECyT-UNC, and CAI+D-UNL. The third author was financed in part by CAPES-Finance Code 001,  CNPq-Project 307261/2023-9, and FAPERJ - CNE 204.037/2024.

\bibliographystyle{siam}

\end{document}